\input amstex
\documentstyle{amsppt}
\topmatter
\magnification=\magstep1
\pagewidth{5.2 in}
\pageheight{6.7 in}
%\hcorrection{-0.4in}
%\vcorrection{-0.4in}
\abovedisplayskip=10pt
\belowdisplayskip=10pt
\NoBlackBoxes
\title
$q$-analogues of the Sums of powers of consecutive integers
\endtitle
%Use \endgraf to indicate new paragraph
%\thanks will become a 1st page footnote
\author  Taekyun Kim \endauthor
\affil\rm{{Institute of Science Education,}\\
{ Kongju National University, Kongju 314-701, S. Korea}\\
{e-mail: tkim$\@$kongju.ac.kr (or tkim64$\@$hanmail.net)}}\\
\endaffil

\abstract{ Let $n, k$ be the positive integers $( k>1 ),$ and let
$S_{n,q}(k)$ be the sums of the $n$-th powers of positive
$q$-integers up to $k-1$: $S_{n,q}(k)=\sum_{l=0}^{k-1}q^l l^n $.
Following an idea due to J. Bernoulli, we explore a formula for
$S_{n,q}(k)$. }
\endabstract
\thanks 2000 Mathematics Subject Classification  11S80, 11B68 \endthanks
\thanks Key words and phrases:  Sums of powers, Bernoulli number, Bernoulli polynomials \endthanks
\rightheadtext{ Taekyun Kim    } \leftheadtext{ Sums of powers of
consecutive $q$-integers  }
%\TagsOnRight
\endtopmatter

\document

\head 1. Introduction \endhead

In the early of 17th century Faulhaber [3, 6, 11] computed the
sums of powers $1^m+2^m+\cdots+n^m$ up to $m=17$ and realized that
for odd $m$, it is not just a polynomial in $n$ but a polynomial
in the triangular number $N=n(n+1)/2$. A good account of
Faulhaber's work was given by Knuth [3, 6, 11]. In 1713, J.
Bernoulli first discovered the method which one can produce those
formulae for the sum $\sum_{l=1}^nl^k$ for any natural numbers
$k$, cf. [1, 2, 4, 5, 13]. Let $q$ be an indeterminate which can
be considered in complex number field, and for any integer $k$
define the $q$-integer as $[k]_q=\frac{q^k-1}{q-1}, $ cf. [2, 7,
8, 9, 10 ]. Note that $\lim_{q\rightarrow 1}[k]_q=k .$ Recently,
many authors are studying  the problems of $q$-analogues of sums
of powers of consecutive integers [3, 4, 5, 6, 7, 8, 11, 12, 13,
14 ]. In this paper we consider the $q$-analogues of the sums of
powers of consecutive integers due to J. Bernoulli .  For any
positive integers $n, k(>1),$ let
$S_{n,q}(k)=\sum_{l=0}^{k-1}q^ll^n.$ Following an idea due to J.
Bernoulli, we explore a formula for $S_{n,q}(k)$ as follows:
$$ q^{-k}S_{l-1,q}(k)+(q^{-k}\log q
)S_{l,q}(k)=\frac{1}{l}\sum_{i=0}^{l-1}
\binom{l}{i}B_{i,q}k^{l-i}+\frac{(1-q^{-k})B_{l,q}(0)}{l},$$ where
$B_{i,q}$ are the $q$-analogue of Bernoulli numbers.

\head 2. $q$-analogue of the Sums of the $n$-th powers of positive
integers up to $k-1$
\endhead

Let $j$ be the positive integers. Then we easily see that
$$q^{j+1}(j+1)-q^jj=(q-1)q^j j +q^{j+1}.\tag 1$$
From Eq.(1), we note that
$$(q-1)\sum_{j=0}^{k-1}q^j j +q[k]_q=q^k k. \tag2$$
Thus, we have the following:
$$S_{1,q}(k)=\sum_{j=0}^{k-1}q^jj=\frac{q^k k-q[k]_q}{q-1}.\tag3$$
Note that $\frac{k(k-1)}{2}=\lim_{q\rightarrow 1}\frac{q^k
k-q[k]_q}{q-1}=S_1(k).$ By the same method of Eq.(1), we easily
see that
$$q^{j+1}(j+1)^2-q^j j^2 =(q-1)q^jj^2+2q^{j+1}j+q^{j+1}.\tag4$$
Thus, we obtain
$$S_{2,q}(k)=\sum_{j=0}^{k-1}q^jj^2=\frac{q^kk^2}{q-1}-2q\frac{q^kk-q[k]_q}{(q-1)^2}-\frac{q[k]_q}{q-1}. \tag5$$
From the simple calculation, we note that
$$\aligned
&q^{j+1}(j+1)^3-q^jj^3=(q-1)q^jj^3+3qq^jj^2+3qq^jj+q^{j+1} ,\\
&q^kk^3=(q-1)S_{3,q}(k)+3qS_{2,q}(k)+3qS_{1,q}(k)+q[k]_q.\endaligned\tag6$$
Thus, we see that
$$S_{3,q}(k)=\frac{q^kk^3}{q-1}-\frac{3q}{q-1}S_{2,q}(k)-\frac{3q}{q-1}S_{1,q}(k)-\frac{q[k]_q}{q-1}.\tag7$$
Let $n, k$ be the positive integers $(k>1)$ and let
$S_{n,q}(k)=\sum_{l=0}^{k-1}q^ll^n. $ Then we have
$$q^{l+1}(l+1)^n-l^{n+1}q^l=q^{l+1}\sum_{i=0}^n\binom{n+1}{i}l^i+(q-1)l^{n+1}q^l .$$
Summing over $l$ from 0 to $k-1$, the left-hand side becomes
$q^kk^{n+1}$, but the right-hand side is a linear combination of
$S_{i,q}(k)$ as follows:
$$q^kk^{n+1}=q(n+1)S_{n,q}(k)+q\sum_{i=0}^{n-1}\binom{n+1}{i}S_{i,q}(k)+(q-1)S_{n+1,q}(k).\tag8$$
Thus, we obtain the following :
\proclaim{ Theorem A} Let $n$, $k
(>1)$ be the positive integers. Then we have
$$S_{n,q}(k)=\frac{k^{n+1}}{n+1}q^{k-1}-\frac{1}{n+1}\sum_{i=0}^{n-1}\binom{n+1}{i}S_{i,q}(k)+\frac{q-1}{q(n+1)}S_{n+1
,q}(k).\tag9$$
\endproclaim

Remark. Note that
$$\aligned
\lim_{q\rightarrow 1}S_{n,q}(k)&=\lim_{q\rightarrow
1}\left(\frac{k^{n+1}}{n+1}q^{k-1}-\frac{1}{n+1}\sum_{i=0}^{n-1}\binom{n+1}{i}S_{i,q}(k)+\frac{q-1}{q(n+1)}S_{n+1
,q}(k)\right)\\
&=\frac{k^{n+1}}{n+1}-\frac{1}{n+1}\sum_{i=0}^{n-1}\binom{n+1}{i}S_{i}(k)=S_n(k),
\text{ cf. [13]}.
\endaligned$$

\head 3. A formula for $S_{n,q}(k)$
\endhead
 In this section, we assume $q\in \Bbb C$ with $|q|<1 .$  Now, we consider the
 $q$-Bernoulli polynomials  as follows:
$$\frac{\log q +t}{qe^t-1}e^{xt}=\sum_{n=0}^{\infty}B_{n,q}(x)\frac{t^n}{n!}, \text{ see [7, 8]. } \tag10$$ In the case $x=0,$
$B_{n,q}(=B_{n,q}(0))$ will be called the $q$-Bernoulli numbers,
see [8, 9]. Let $F_q(t,x)=\frac{\log q +t}{qe^t-1}e^{xt}.$ Then we
see that
$$F_q(t,x)=-(t+\log q)\sum_{n=0}^{\infty}e^{(n+x)t}q^n, \text{
$t\in\Bbb C$ with $|t|<2\pi$, cf.[8, 9, 10 ]. }\tag11$$ In [8], it
was known that
$$B_{n,q}(x)=\sum_{j=0}^{n}\binom{n}{j}B_{j,q}x^{n-j}=m^{n-1}\sum_{i=0}^{m-1}q^iB_{n,q^m}(\frac{x+i}{m}), \text{ for
$n\geq 0$}.\tag12$$
 By (10) and (11), we easily see that  the
$q$-Bernoulli numbers can be rewritten as
$$B_{0,q}=\frac{q-1}{\log q}, \text{ }
q(B_q +1)^k-B_{k,q}=\delta_{k,1}, \text{ for $k\geq 1 , $}
$$ where $\delta_{k,1} $ is Kronecker symbol and we use the usual
convention about replacing $B_q^{i}$ by $B_{i,q}$, cf. [8, 9, 10].
From the simple calculation, we note that
$$-\sum_{n=0}^{\infty}e^{(n+k)t}q^n+\sum_{n=0}^{\infty}e^{nt}q^{n-k}
=\sum_{l=0}^{\infty}(q^{-k}\sum_{n=0}^{k-1}n^lq^n)\frac{t^l}{l!}
=\sum_{l=1}^{\infty}(q^{-k}l\sum_{n=0}^{k-1}n^{l-1}q^n)\frac{t^{l-1}}{l!}.\tag13$$
Thus, we have
$$\aligned
 \sum_{l=0}^{\infty}(B_{l,q}(k)-q^{-k}B_{k,q}(0))\frac{t^l}{l!}
&=-(t+\log q )\sum_{n=0}^{\infty}e^{(n+k)t}q^n +(t+\log
q)\sum_{n=0}^{\infty}e^{nt}q^{n-k}\\
&=\sum_{l=0}^{\infty}\left(q^{-k}lS_{l-1,q}(k)+(q^{-k}\log
q)S_{l.q}(k)\right)\frac{t^l}{l!}.
\endaligned\tag14$$
By comparing the coefficient on both side in Eq.(14), we obtain
the following:

\proclaim{ Theorem B} Let $l,k$ be positive integers ( $k>1$).
Then we see that
$$q^{-k}S_{l-1,q}(k)+(q^{-k}\log
q )S_{l,q}(k)=\frac{B_{l,q}(k)-q^{-k}B_{l,q}(0)}{l}.$$
\endproclaim
Remark. Note that
$$S_{l-1}(k)=\lim_{q\rightarrow 1}\left(q^{-k}S_{l-1,q}(k)+(q^{-k}\log
q )S_{l,q}(k)\right)=\frac{B_l(k)-B_l(0)}{l},\text{ cf. [13],}$$
where $B_l(k)$ are called ordinary Bernoulli polynomials.

\enddocument